\documentclass[11pt]{amsart}

\usepackage{amsmath}
\usepackage{amsfonts}
\usepackage{amssymb}
\usepackage{amsthm}
\usepackage{euscript}
\usepackage{verbatim}
 \usepackage{epsfig}
\newcommand{\Hmm}[1]{\leavevmode{\marginpar{\tiny%
$\hbox to 0mm{\hspace*{-0.5mm}$\leftarrow$\hss}%
\vcenter{\vrule depth 0.1mm height 0.1mm width \the\marginparwidth}%
\hbox to 0mm{\hss$\rightarrow$\hspace*{-0.5mm}}$\\\relax\raggedright #1}}}

\newcommand{\Rr}{{\bf R}}
\newcommand{\reals}{{\bf R}}

\newcommand{\nrm}{\Arrowvert}
\newcommand{\kuc}{\lesssim}

\newcommand{\buy}{\gtrsim}

\makeatletter

\def\e{\varepsilon}
\def\eps{\varepsilon}

\def\scriptc{{\mathcal C}}
\def\scriptf{{\mathcal F}}

\def\rp{{p^{-1}}}
\def\rq{{q^{-1}}}
\def\rr{{r^{-1}}}

\newtheorem{theorem}{Theorem}
\newtheorem{lemma}{Lemma}
\newtheorem{prop}[lemma]{Proposition}

\newtheorem{definition}{Definition}

\newcounter{smalllist}
\renewcommand{\thesmalllist}{{\rm\roman{smalllist}}}

\title[Mixed norm estimates]{Mixed norm estimates for certain generalized Radon transforms}

\author{Michael Christ} 
 \address{Department of Mathematics\\
 University of California\\
 Berkeley, CA 94720-3840}

\author {M. Burak Erdo\smash{\u{g}}an}
\address{Department of Mathematics \\
University of Illinois \\
Urbana, IL 61801}

\date{9/6/2005}
\thanks{The second author was partially supported by the NSF grant DMS-0540084.}

\begin{document}


\maketitle

\section{Introduction}

In this paper we investigate the mapping properties in Lebesgue-type spaces
of certain generalized Radon transforms defined by integration over curves. 

Let $X$ and $Y$ be open subsets of $\Rr^d$, $d\geq 2$, and let $Z$ be a smooth submanifold of 
$X\times Y\subset  \Rr^{2d}$ 
of dimension $d + 1$. Assume that 
the projections $\pi_1:Z\rightarrow X$ and $\pi_2:Z\rightarrow Y$ are submersions at each point of $Z$.
For each $y\in Y$, let 
$$
\gamma_y=\{x\in X: (x,y)\in Z\}=\pi_1\pi_2^{-1}(y).
$$ 
In this case, $\gamma_y$ are smooth curves in $X$
which vary smoothly with $y\in Y$.
For every $y\in Y$, choose a smooth, non-negative measure $\sigma_y$ on $\gamma_y$ 
which varies smoothly with $y$ in the natural sense.
A generalized Radon transform $T$ 
(see e.g.\ \cite{christ2, greenleafseegerwainger, seeger})
is defined as an operator taking functions on $X$ to functions on $Y$ via
$$
Tf(y)=\int_{\gamma_y} f d\sigma_y.
$$ 
The adjoint of this operator has a similar form:
$$
T^*g(x)=\int_{\gamma^*_x} g d \sigma^*_x,
$$
where
$$
\gamma^*_x=\{y:(x,y)\in Z\}=\pi_2\pi_1^{-1}(\{x\})\subset Y
$$
and $\sigma_x^*$ is a nonnegative measure on $\gamma_x^*$ with a smooth density
which varies smoothly with $x$.

Tao and Wright \cite{taowright} have formulated and proved a nearly optimal 
characterization of the local $(L^p,L^q)$ 
mapping properties of these operators. We extend their result to the mixed-norm 
setting and  
obtain essentially optimal local mixed-norm inequalities for  
these operators, under one additional dimensional restriction.
Previously this result was obtained  
for a model operator in \cite{wolff3d,erdogan,christerdogan}.  
See \cite{christcomb,christ3,christerdogan, erdogan, taowright,wolff3d} 
for various examples and prior work.

\paragraph{\bf Mixed norms on $Y$.}  

Throughout the entire discussion, $X,Y$ will denote sufficiently small neighborhoods
of $x_0,y_0$ for some fixed point $(x_0,y_0)\in Z$.
Let $\Pi:Y\rightarrow \Rr$ be a submersion such that the fibers $\Pi^{-1}(t)$ are transverse to the 
curves $\gamma^*_x$. This means that the restriction of $\Pi$ to
$\gamma^*_x$ is a diffeomorphism for each $x\in X$.
Choose coordinates so that $\Pi(y_0)=0$.
Let $\lambda_t$ be Lebesgue measure on the $d-1$-dimensional surface $\Pi^{-1}(t)$. 
To a function $f:Y\mapsto \Rr$ we associate the mixed norms
$$
\|f\|_{L^qL^r(Y)}=
\nrm f\nrm_{q,r}
:=\left[\int_{\Rr}\left[\int_{\Pi^{-1}(t)} |f(s)|^r d \lambda_t(s)\right]^{q/r}dt\right]^{1/q}.
$$
The integral with respect to $t$ is taken over a small neighborhood of the origin
in $\reals$; we may assume that this neighborhood is contained in $[-1,1]$.
Here, and throughout this paper, $t$ is restricted to lie in a one-dimensional manifold.
This is not a natural restriction, but our analysis yields reasonably satisfactory
results only in that special case. 

We say that $T$ is of strong mixed type $(p,q,r)$ if $T$ maps $L^p(X)$ to $L^qL^r(Y)$ boundedly.
We are mainly interested in
local estimates. 
We assume throughout the discussion  that $T$ is $L^p$-improving, 
which means that for each $p\in (1,\infty)$,
there exists $q>p$ such that $T$ maps $L^p(X)$ to $L^q(Y)$. See \cite{cnsw} and \cite{seeger} 
for characterizations of this property.

Our theorem is a characterization of the exponents $(p,q,r)$ for which $T$ is bounded.
Before we proceed let us record the simple facts about these exponents. \\
{\bf (i)} Because of the transversality hypothesis described above, $T$ is   of strong mixed type 
$(p,\infty,p)$ for all $p\in [1,\infty]$.\\
{\bf (ii)} Since we are working in a bounded region, whenever $T$ is of strong mixed type
$(p,q,r)$, it is also of strong mixed type
$(p_1,q_1,r_1)$ whenever $p_1\geq p$, $q_1\leq q$, and $r_1\leq r$.\\
{\bf (iii)} Because of  (i) and (ii), 
$T$ is of strong mixed type $(p,q,r)$ whenever $p\geq r$.

\paragraph{\bf Two-parameter Carnot-Carath\'{e}odory Balls.} 
 
Tao and Wright 
\cite{taowright} 
related the set of all 
exponents $(p,q)$ for which
the operator $T$ maps   $L^p(X)$ to $L^q(Y)$ to 
the geometry of $Z$. To describe this relation, 
choose smooth nowhere-vanishing linearly independent real vector fields $V_1$, $V_2$ on $Z$ 
whose integral curves are the fibers of 
$\pi_1$, $\pi_2$ respectively. Equivalently, at each point $z\in Z$, $V_j$ spans the 
nullspace of $D\pi_j$, for $j=1,2$.
The $L^p$-improving property is equivalent to
$V_1,V_2$ satisfying the bracket condition, i.e., $V_1$ and $V_2$ together with their iterated 
commutators span the tangent space to $Z$ at each point in $Z$ \cite{cnsw}.

\begin{definition}\label{D:balls}
 Let $z_0\in Z$ and $0<\delta_1, \delta_2\ll 1$. The {\it two-parameter Carnot-Carath\'{e}odory ball}
$B(z_0,\delta_1,\delta_2)$ consists of all the points $z\in Z$ such that there exists  an 
absolutely continuous function $\varphi:[0,1]\rightarrow Z$ satisfying\\
(i) $\varphi(0)=z_0$, $\varphi(1)=z$\\
(ii) for almost every $t\in [0,1]$
$$
\varphi^\prime(t)=a_1(t)V_1(\varphi(t))+a_2(t)V_2(\varphi(t))
$$ 
with $|a_1(t)|<\delta_1$, $|a_2(t)|<\delta_2$. 
\end{definition} 

The metric properties of 
Carnot-Carath\'{e}odory balls were studied extensively in \cite{nagelsteinwainger}.
The discussion there is phrased in terms of the one-parameter family of balls naturally associated
to a family of vector fields satisfying the bracket condition. These balls depend on a center
point, a radius $r$, and a family $\{W_j\}$ of vector fields. 
They can equivalently be viewed as depending on a center point and the family $\{rW_j\}$ of vector fields,
with the radius redefined to be identically one. In these terms, the proofs in 
\cite{nagelsteinwainger} go through more generally, for balls of radius one assigned to
families of vector fields $\{U^\alpha_j: 1\le j\le J\}$ satisfying the bracket condition for each
parameter $\alpha$, with
appropriate uniformity as $\alpha$ varies.
In particular, for the vector fields $\{\delta_1V_1,\delta_2V_2\}$,
provided that $0<\delta_1,\delta_2\le c_0$ for sufficiently small $c_0$,
the conclusions of \cite{nagelsteinwainger} hold uniformly in $\delta_1,\delta_2$ under
a supplementary hypothesis of weak comparability, which is discussed below.
See \cite{taowright, christ3}. 

It will be convenient in our proof to
parametrize the curves $\gamma^*_x$ by $t$ so that
$\Pi(\gamma^*_x(t))\equiv t$.
With this parametrization, the measure $\sigma_x$ on $\gamma^*_x$ is equivalent to $dt$,
uniformly in $x$.
We rescale $V_1$ if necessary so that for each $z \in Z$ and sufficiently small $s\in \Rr$ 
\begin{equation}\label{v1flow}
\Pi\pi_2(e^{sV_1}z)=\Pi\pi_2(z)+s.
\end{equation}

\begin{definition}
Let $0<\theta\le 1$, and let $A$ be a positive constant.
We say that $0<\delta_1, \delta_2\ll 1$ are $(\theta,A)$-weakly comparable,
and write $\delta_1 \sim_{\theta,A}  \delta_2$,
if $\delta_1 \leq A\delta_2^\theta$ and $\delta_2 \leq A\delta_1^\theta$. 
\end{definition}
The following lemma collects
basic facts about the balls  $B(z,\delta_1,\delta_2)$.
 
\begin{lemma}\label{L:balls}
Let $K$ be a compact subset of $Z$. Assume that 
$\delta_1 \sim_{(\theta,A)} \delta_2$ are sufficiently small,
and let $z\in K$.   
Then $B=B(z,\delta_1,\delta_2)$ satisfies\\
(i) $|B|\sim|B(z,2\delta_1,2\delta_2)|$,\\
(ii) $|B|\sim|\pi_1(B)|\delta_1\sim|\pi_2(B)|\delta_2$\\
(iii) $|\Pi(\pi_2(B))|\sim\delta_1$,\\
(iv) $\nrm \chi_{\pi_2(B)}\nrm_{q^\prime,r^\prime}\sim
|B|^{1-\frac{1}{r}}\delta_1^{\frac{1}{r}-\frac{1}{q} }\delta_2^{\frac{1}{r}-1}$,\\
(v) 
$$
\frac{|B|}{|\pi_1(B)|^{\frac{1}{p}} \nrm\chi_{\pi_2(B)}\nrm_{q^\prime,r^\prime}}\sim
|B|^{\frac{1}{r}-\frac{1}{p}}\delta_1^{\frac{1}{p}+\frac{1}{q}-\frac{1}{r}}\delta_2^{1-\frac{1}{r}}.
$$
Here $1\leq p,q,r\leq\infty$, and $q^\prime, r^\prime$ are the exponents conjugate to
$q,r$, respectively. 
\end{lemma}
\noindent The notation $A\sim C$ means that the ratio $A/C$ is bounded above and
below by quantities depending 
on $Z$, $\theta$, $A$, and the compact set $K$, but not on $\delta_1,\delta_2$.
In the absence of weak comparability,
the doubling property (i) fails in general for two-parameter Carnot-Carath\'eodory
balls associated to $C^\infty$ vector fields satisfying the bracket condition
\cite{christ3}.

For a sketch of the proof of the lemma see \S\ref{formerlyAppendix1} below.

\paragraph{\bf Statement of results.} 

Recall that $T$ is said to be of restricted weak type $(p,q)$ if for all Lebesgue measurable sets
$E\subset X$ and $F\subset Y$
\begin{equation}\label{nonmixed}
\langle T\chi_E,\chi_F\rangle \kuc |E|^{1/p} |F|^{1/q^{\prime}},
\end{equation}
where $q^\prime$ denotes the exponent conjugate to $q$.
In our setup 
\begin{equation}
\label{pairing}
\langle T\chi_E,\chi_F\rangle \approx |\pi_1^{-1}(E)\cap\pi_2^{-1}(F)|,
\end{equation}
where $|\cdot|$ denotes Lebesgue measure on $Z$. We test the inequality 
\eqref{nonmixed} on the Carnot-Carath\'{e}odory balls $B(z,\delta_1,\delta_2)$
under the restriction that $\delta_1\sim_{(\theta,A)}\delta_2$. 
Let $E=\pi_1(B(z,\delta_1,\delta_2)), F=\pi_2(B(z,\delta_1,\delta_2))$. 
Using \eqref{pairing}, Lemma~\ref{L:balls} and restricting attention 
to the nontrivial case where $q>p$, 
the inequality \eqref{nonmixed} reads   
\begin{equation}\label{c1c2}
|B(z,\delta_1,\delta_2)|\buy \delta_1^{c_1}\delta_2^{c_2},
\end{equation}
where 
\begin{equation}\label{c1c2pq}
c_1=\frac{\rp}{\rp-\rq },\qquad
c_2=\frac{1-\rq}{\rp-\rq}.
\end{equation}

Define
$$
\scriptc_{\theta,A}(T):=\{(c_1,c_2): 
\inf\left(\frac{|B(z,\delta_1,\delta_2)|}{\delta_1^{c_1}\delta_2^{c_2}}\right)>0 
 \},
$$
where the infemum is taken over all $z\in Z$ and over all pairs $\delta_1,\delta_2$
that satisfy $\delta_1\sim_{(\theta,A)}\delta_2$. 
Define
$$
\scriptc(T):=\cap_{0<\theta\le 1} \cap_{A\ge 1}\scriptc_{\theta,A}(T).
$$

According to \eqref{c1c2}, \eqref{nonmixed} 
can not hold for $(p,q)$ if the corresponding $(c_1,c_2)$
does not belong to  $\scriptc(T)$. 
Tao and Wright \cite{taowright} proved 
that\footnote{Tao and Wright defined the set $\scriptc(T)$ differently. 
An analysis of the two-parameter balls along the lines of \cite{nagelsteinwainger} 
establishes the equivalence of these two definitions.}
for all $(c_1,c_2)$ in the interior of $\scriptc(T)$, 
\eqref{nonmixed} holds for the exponents $(p,q)$ defined by \eqref{c1c2}.

In this note we  extend this result to mixed norms. 
We say that $T$ is of restricted weak mixed type $(p,q,r)$ if
for all $E\subset X$ and $F\subset Y$,
$$
\langle T\chi_E,\chi_F\rangle \kuc |E|^{1/p} \nrm \chi_F\nrm_{q^\prime,r^\prime}.
$$
By interpolation, the strong mixed type estimates can be obtained from  these inequalities,
except for exponents corresponding to boundary points of $\scriptc(T)$.

The two-parameter Carnot-Carath\'{e}odory balls defined above 
also dictate the allowed exponent triples $(p,q,r)$ for mixed norm inequalities,
under certain additional restrictions on the exponents $p,q,r$: 
\begin{definition} \label{D:ptheta}
Let $P_{\theta,A}(T)$ be the set of all exponents  $(p,q,r)$ satisfying \\
(i) $1\leq p\leq q\leq r\leq \infty$,  \\
(ii) 
\begin{equation}
\sup_{z,\delta_1,\delta_2} 
\frac{|B(z,\delta_1,\delta_2)|}{|\pi_1(B(z,\delta_1,\delta_2))|^{1/p}\nrm \chi_{\pi_2(B(z,\delta_1,\delta_2))}\nrm_{q^\prime,r^\prime}}<\infty,
\end{equation}
where
$q^\prime,r^\prime$ are the conjugates of $q,r$ respectively, 
and the supremum is taken over all $z\in Z$ and 
$\delta_1 \sim_{\theta,A}\delta_2$. 
\end{definition}

Using Lemma~\ref{L:balls}
we can rewrite the second condition in the definition of $P_{\theta,A}(T)$ as in \eqref{c1c2} with
\begin{equation}\label{c_1c_2pqr}
c_1=\frac{\rp+\rq-\rr}{\rp-\rr },\qquad
c_2=\frac{1-\rr}{\rp-\rr}.
\end{equation}
Define
$$
P(T):=\cap_{0<\theta\le 1}\cap_{A\ge 1}
P_{\theta,A}(T) =\{(p,q,r):r\geq q \geq p, (c_1,c_2)\in \scriptc(T)\}.
$$
We assume always that $r>p$, since otherwise the desired inequality
holds automatically so long as $\pi_1,\pi_2$ are submersions, as discussed above.

It is natural to conjecture that if $1\leq p\leq q \leq r \leq \infty$, 
then $T$ is of restricted weak mixed
type $(p,q,r)$ if and only if $(p,q,r)\in P(T)$; 
we prove this conjecture except at the endpoints.
In an appendix we explain the presence of the
additional restriction $p\leq q\leq r$. 

\begin{theorem} \label{T:main} Let $T$ be a generalized Radon transform
of the class described above,  
and let $1\le p\le q\le r\le\infty$.
If $(p,q,r)$ is in the interior of $P(T)$, then $T$ maps  $L^p(X)$ to $L^qL^r(Y)$ boundedly.
Moreover, if $(p,q,r)\not\in P(T)$,
then $T$ does not map $L^p(X)$ to $L^qL^r(Y)$ boundedly.
\end{theorem}

\section{A brief account of \cite{taowright}.} 

Given $E\subset X$ and , $F\subset Y$, let 
\begin{equation}\label{alphabeta}
\alpha_1=\frac{\langle T\chi_E,\chi_F\rangle}{|E|},\,\,\,\,\,
\alpha_2=\frac{\langle T\chi_E,\chi_F\rangle}{|F|}.
\end{equation}
\eqref{nonmixed} follows from an estimate of the form
$$
|\Omega|:=
|\pi_1^{-1}(E)\cap\pi_2^{-1}(F)|
\approx  \langle T\chi_E,\chi_F\rangle \buy \alpha_1^{c_1}\alpha_2^{c_2},
$$
where $c_1,c_2$ are as in \eqref{c1c2pq}. Moreover, a loss of an arbitrarily small 
power of $\alpha_1\alpha_2$ is of no consequence
since $(c_1,c_2)$ belongs to the {\em interior} of $\scriptc(T)$. 
Throughout the remainder of the discussion we will denote
\begin{equation}
\Omega:= \pi_1^{-1}(E)\cap\pi_2^{-1}(F).
\end{equation}

Let $\e>0$ be a small exponent to be specified near the end of the proof,
and let $C_\eps$ be a large constant. 
All constants in the ensuing discussion
depend on $X,Y,Z,c_1,c_2$ and the quantities given as subscripts. 
\begin{definition}
Let $S\subset[-1,1]$ be a measurable set of positive Lebesgue measure.
We say that $S$  is central with width $w>0$ if \\
(i)  $S\subset [-C_\eps w,C_\eps w]$ and\\
(ii) $|I\cap S|\leq C_{\eps} (|I|/w)^\eps |S|$ for all intervals $I$.
\end{definition}

Let $\Omega$ be a Lebesgue measurable subset of $Z$, having positive measure. 
We will assume throughout the discussion that $\pi_j(\Omega)$ has positive measure,
for both $j=1,2$.
Define
\begin{gather}
\label{alphasdefn}
\alpha_j=\frac{|\Omega|}{|\pi_j(\Omega)|} \text{ for $j=1,2$} 
\\
\alpha:=\min(\alpha_1,\alpha_2).
\end{gather}

Fix any $\rho>0$; $\rho$ will eventually  be chosen to be arbitrarily small
at the conclusion of the proof.
\begin{definition}
Let $\Omega\subset Z$ and let $\alpha_1,\alpha_2$ be as defined in \eqref{alphasdefn}.
Let $j$ be an integer. A {\it $j$-subsheaf} $\Omega^\prime$ of $\Omega$,
of width $w_j$, is a subset of $\Omega$ of measure 
$\geq C_{\rho}^{-1} \alpha^{\rho} |\Omega|$ such that for all $x\in \Omega^\prime$, the set
$$
\{|t| \ll 1:e^{tV_j}(x)\in \Omega^\prime \}
$$
is a central set of width $w_j$ and measure $\geq C_{N,\e} \alpha^{ C_N\e+C/N} \alpha_j$.
Here $V_j=V_1$ of $j$ is odd, and $V_j=V_2$ if $j$ is even.
\end{definition}

In \cite{taowright} it is  proved (see Corollary 8.3) that
for any $\rho\in(0,1]$ there exists $C_\rho<\infty$
such that
for any measurable set $\Omega\subset Z$ there exists a nested sequence of subsets of $\Omega $ 
$$
\Omega_0\subset\Omega_1\subset...\subset\Omega_{d+1}\subset \Omega 
$$
such that:
For each $j$, $\Omega_j$ is a $j$-sheaf of $\Omega $ with width   $w_j$,
$$
C_{\rho}^{-1}\alpha^{\rho}\delta_j\leq w_j \leq \delta_j,
$$
where $\delta_j$ is a 2-periodic sequence (that is, $\delta_{j+2}=\delta_j$
for all $0\le j\le d-1$) with the properties
$$
C_{\rho}^{-1}\alpha^{\rho}\alpha_{ j}\leq \delta_j \leq 1,
$$
and
$$
\delta_1\le C_\rho\delta_2^\rho,
\qquad
\delta_2\le C_\rho\delta_1^\rho.
$$
In particular, $\delta_1,\delta_2$
are weakly comparable, even though $\alpha_1,\alpha_2$ need not be; 
this ultimately explains why only balls with weakly comparable radii $\delta_1,\delta_2$
need be taken into account in the hypothesis of our theorem.

Using this construction, it is proved \cite{taowright} that there exists some ball 
$B(z,\delta_1,\delta_2)$  
such that   
\begin{equation}
|B(z,\delta_1,\delta_2)\cap \Omega| 
\ge c
\alpha^{\varrho}
\left(\frac{\alpha_{1}}{\delta_1}\right)^{\lfloor(d+2)/2\rfloor}
\left(\frac{\alpha_{2}}{\delta_2}\right)^{\lfloor(d+1)/2\rfloor}
|B(z,\delta_1,\delta_2)|, 
\end{equation}
where $\varrho>0$ can be made arbitrarily small by choosing $\rho$ sufficiently small,
and where $c>0$ depends on $\varrho$ but not on $\Omega,\alpha_j,\delta_j$.
This implies that for arbitrarily small $\varrho>0$ 
\begin{equation}\label{omegabound1}
|\Omega|\buy C_{\varrho}\alpha^{\varrho} \delta_1^{c_1}\delta_2^{c_2}
\left(\frac{\alpha_{1}}{\delta_1}\right)^{\lfloor(d+2)/2\rfloor}
\left(\frac{\alpha_{2}}{\delta_2}\right)^{\lfloor(d+1)/2\rfloor},
\end{equation}
for all $(c_1,c_2)$ in the interior of $\scriptc(T)$.
This finishes the proof in the non mixed-norm case since\footnote{
The bounds $c_1,c_2\ge d$ are a simple consequence of an 
equivalent definition of $\scriptc(T)$ given in \cite{taowright}.}
$c_1,c_2\geq d \geq \lfloor(d+2)/2\rfloor, \lfloor(d+1)/2\rfloor.$
The roles of $X$ and $Y$ are interchangeable. 
Therefore there is the alternative bound
\begin{equation}\label{omegabound}
|\Omega|\buy C_{\varrho}\alpha^{\varrho} \delta_1^{c_1}\delta_2^{c_2}
\left(\frac{\alpha_{1}}{\delta_1}\right)^{\lfloor(d+1)/2\rfloor}
\left(\frac{\alpha_{2}}{\delta_2}\right)^{\lfloor(d+2)/2\rfloor},
\end{equation}
for all $(c_1,c_2)$ in the interior of $\scriptc(T)$.

\section{Proof of Theorem~\ref{T:main}.}

The second conclusion of the theorem
is immediate. If $(p,q,r)\not \in P(T)$, then $(p,q,r)\not \in P_{\theta,A}(T)$ 
for some $\theta,A$.  Therefore $T$ can not be of mixed type $(p,q,r)$.

For the first conclusion, it suffices to prove that
for all $(p,q,r)$ in the interior of $P(T)$ and for all $\eta>0$,  $\beta>0$,
and  $F\subset Y$,
\begin{equation}\label{main}
\langle T\chi_E,\chi_F\rangle\leq C_{p,q,r,\eta} 
\beta^{-\eta}|E|^{1/p}\nrm \chi_F\nrm_{q^\prime,r^\prime},
\end{equation}
where  $E:=\{x\in X:  \beta<T^*\chi_F(x)\leq 2\beta\}$.

Let $\Omega:=\pi_1^{-1}(E)\cap \pi_2^{-1}(F)$. 
Note that 
\begin{equation} \label{omega}
|\Omega| \approx\langle T\chi_E,\chi_F\rangle 
=\langle \chi_E, T^*\chi_F \rangle \approx \beta |E|.
\end{equation}
For each $x\in E$, let 
$\scriptf(x):=\{t\in\Rr: \gamma^*_x(t) \in F\}\subset \Pi(F)
\subset [-1,1]$. Note that $|\scriptf(x)|\approx \beta$.  

Let $\eta>0$ be a small constant to be specified at the end of the proof.
Fix a small constant $c_\eta>0$. 
Let $I(x)\subset [-1,1]$ be a dyadic interval of minimal length so that 
\begin{equation}
\label{I(x)defn}
|I(x)\cap \scriptf(x)|\geq c_\eta |I(x)|^\eta |\scriptf(x)|;
\end{equation}
we choose $c_\eta$ to guarantee that 
either $I=[-1,0]$ or $I=[0,1]$ satisfies \eqref{I(x)defn}. 
No interval of length $<(c_\eta|\scriptf(x)|)^{1/(1-\eta)}$ can satisfy \eqref{I(x)defn},
so there must exist at least one dyadic interval of minimal length 
among all those satisfying the inequality.

This implies that for any dyadic subinterval $J\subset I(x)$,
\begin{equation}\label{localization}
|J\cap \scriptf(x)|\leq \frac{|J|^\eta}{|I|^\eta}|I(x)\cap \scriptf(x)|.
\end{equation}
Let 
$$
E_{m,k}=\{x\in E:|I(x)|\approx 2^m\beta, |I(x)\cap \scriptf(x)|\approx 2^k\}.
$$
Note that $\beta^{\eta/(1-\eta)}\kuc 2^m \kuc \beta^{-1}$ 
and $\beta^{1/(1-\eta)}\kuc 2^k\kuc \beta$.

By the pigeonhole principle
there exists a pair $m,k$ 
such that $\tilde E = E_{m,k}$ satisfies
\begin{equation}
\langle T\chi_E,\chi_F\rangle
\le C_\eta\beta^{-\eta}
\langle T\chi_{\tilde E},\chi_F\rangle.
\end{equation}
Choose one such pair, with which we work exclusively henceforth. 


Partition $\Rr$ into 
intervals $I_n$, each of length $C2^m\beta$
where $C$ is a sufficiently large fixed constant.
Set 
$$
E_n=\{x\in \tilde E:I(x)\cap  I_n\neq \emptyset\},
$$
$$
F_n=F\cap \Pi^{-1}(I_{n-1}\cup I_n\cup I_{n+1}).
$$
Note that 
\begin{equation}\label{localized}
\langle T\chi_{\tilde E}, \chi_F\rangle \kuc \sum_n \langle T\chi_{E_n}, \chi_{F_n} \rangle.
\end{equation} 
Let $\Omega^n:=\pi_1^{-1}(E_n)\cap\pi_2^{-1}(F_n)$. 
Then
\begin{equation}\label{omega^n}
2^k|E_n|\lesssim |\Omega^n|\le \beta |E_n|
\end{equation}
by Fubini's theorem. For on one hand, $|\pi_1^{-1}(x)|\le\beta$ for all $x\in E$.
On the other hand, $2^k\sim |I(x)\cap\scriptf(x)|$ and
for $x\in E_n$, 
$I(x)\cap\scriptf(x)\subset (I_{n-1}\cup I_n\cup I_{n+1})\cap\scriptf(x)
\subset \pi_2^{-1}(F_n)$.
Using \eqref{localization} with $I=I_{n-1}\cup I_n \cup I_{n+1}$,   
for each $x\in E_n$ and any dyadic interval $J$, we obtain
\begin{equation}\label{widthbound}
|J\cap \scriptf(x)|\kuc |J|^\eta(2^m\beta)^{-\eta} 2^k.
\end{equation}
Let $\alpha_{n,i}=|\Omega^n|/|\pi_i(\Omega^n)|$, $i=1,2.$
Then 
\begin{equation}
2^k\lesssim \alpha_{n,1}\lesssim \beta,
\end{equation} 
by \eqref{omega^n}. 
Let $\alpha_n=\min(\alpha_{n,1},\alpha_{n,2})$.

The following lemma follows from a simple application of H\"{o}lder's inequality; 
see \cite{christerdogan}.
\begin{lemma}\label{L:holder}
Let $F\subset Y$. For $r\geq q$, 
$$
\nrm\chi_F\nrm_{q^\prime r^\prime}\geq |F|^{1/r^\prime}|\Pi(F)|^{1/q^\prime-1/r^\prime}.
$$
\end{lemma}

We aim to prove that
$
\langle T\chi_{E_n},\chi_{F_n}\rangle \kuc |E_n|^{1/p}\nrm \chi_{F_n}\nrm_{q^\prime,r^\prime}
$.
Via the preceding lemma, this would follow from
\begin{equation}\label{main1}
|\Omega^n| \approx\langle T\chi_{E_n},\chi_{F_n}\rangle\buy
\alpha_{n,1}^{\gamma_1 }\alpha_{n,2}^{\gamma_2 }|\Pi(F_n)|^{\gamma_3},
\end{equation}
where  
\begin{equation} \label{gammasdefn}
\gamma_1=\frac{\rp}{\rp-\rr},\qquad
\gamma_2=\frac{1-\rr}{\rp-\rr},\qquad
\gamma_3=\frac{\rq-\rr}{\rp-\rr}.
\end{equation}
Whereas the exponent $\frac1{q'}-\frac1{r'}$ in Lemma~\ref{L:holder} was nonpositive,
here all the exponents $\gamma_j$ are nonnegative.
 
The construction of Tao and Wright \cite{taowright} produces a nested sequence of subsets 
$ \Omega^n_0\subset\Omega^n_1\subset...\subset\Omega^n_{d+1}\subset \Omega^n $
satisfying:
For each $j$, $\Omega^n_j$ is a $j$-sheaf of $\Omega^n$ with width $w_j$,
$$
C_{\rho}^{-1}\alpha_{n}^{\rho}\delta_j\leq w_j \leq C_{\rho}\delta_j,
$$
where $\delta_j$ is a 2-periodic sequence with the property
\begin{equation}\label{deltaalpha}
C_{\rho}^{-1}\alpha_{n}^{\rho}\alpha_{n,j}\leq \delta_j \leq  1,
\end{equation}
and
$|\Omega^n_j|\ge c_\rho \alpha_{n}^\rho |\Omega^n|.$ 
Moreover, $\delta_1,\delta_2$ are weakly comparable. 

Now, we prove  that the construction above   guarantees a lower bound for $\delta_1$. Indeed, let $j$ be odd.
Using  \eqref{widthbound}, we have
\begin{equation}\label{gecici}
|\Omega^n_j|\leq |E_n| \sup_{J,x} |J\cap F_n(x)|\leq |E_n| w_j^\eta(2^m\beta)^{-\eta} 2^k,
\end{equation}
where the supremum is taken over all intervals   $J$ of length $w_j$ and all points $x\in E_n$. 
Using \eqref{gecici}, \eqref{omega^n} and the fact that $\Omega_j^n$ is a $j$-sheaf of $\Omega^n$, we get
\begin{equation}
c_\rho\alpha_n^{\rho}  2^k |E_n|
\lesssim
C_{\rho} \alpha_n^{\rho} |\Omega^n|
\leq |\Omega^n_j|
\leq |E_n| w_j^\eta(2^m\beta)^{-\eta} 2^k.
\end{equation}
This implies that for odd $j$,
\begin{equation}\label{delta1bound} 
\delta_1\buy  w_j\buy 2^m\beta C_{\rho}^{1/\eta} \alpha_n^{\rho/\eta}
\gtrsim |\Pi(F_n)|C_{\rho}^{1/\eta} \alpha_n^{\rho/\eta}. 
\end{equation}

Using the fundamental information \eqref{omegabound}, we obtain
$$
|\Omega^n|\buy 
\alpha_n^\varrho
\alpha_{n}^{C\rho/\eta}\delta_1^{c_1}\delta_2^{c_2}
\left(\frac{\alpha_{n,1}}{\delta_1}\right)^{\lfloor(d+1)/2\rfloor}
\left(\frac{\alpha_{n,2}}{\delta_2}\right)^{\lfloor(d+2)/2\rfloor},
$$
for any $(c_1,c_2)$ belonging to the interior of $\scriptc(T)$.
Note that for $d\geq 2$ and for such $(c_1,c_2)$, we have
$0\le \gamma_3<1\le c_1-\lfloor (d+1)/2\rfloor$ and   $c_2>\lfloor (d+2)/2\rfloor$. 
Using this, \eqref{deltaalpha}
and \eqref{delta1bound}, we obtain
$$
|\Omega^n|\buy \alpha_n^\varrho \alpha_{n}^{C\rho/\eta}\delta_1^{\gamma_3} 
 \alpha_{n,1}^{c_1-\gamma_3} 
 \alpha_{n,2}^{c_2}  
 \buy \alpha_n^\varrho
\alpha_{n}^{C\rho/\eta}|\Pi(F_n)|^{\gamma_3} 
 \alpha_{n,1}^{c_1-\gamma_3} 
 \alpha_{n,2}^{c_2}, 
$$
for all $(c_1,c_2)$ in the interior of  $\scriptc(T)$. 

Given exponents $(p,q,r)\in P(T)$,
define $\gamma_1,\gamma_2,\gamma_3$ by \eqref{gammasdefn}.
Choose $\rho>0$, depending on $\eta$, and $(c_1,c_2)$
in the interior of $\scriptc(T)$, so that
$$\frac{C\rho}{\eta}+c_1-\gamma_3+\varrho\leq \gamma_1,
\qquad
\frac{C\rho}{\eta}+c_2 +\varrho \leq \gamma_2.$$ 
Therefore,
since $\gamma_3\ge 0$
(which is a consequence of the assumptions that $r>p$ and $r\ge q$), 
$$
|\Omega^n|\buy   \alpha_{n,1}^{\gamma_1}\alpha_{n,2}^{\gamma_2} |\Pi(F_n)|^{\gamma_3},
$$
and hence
$$
\langle T\chi_{E_n}, \chi_{F_n} \rangle \kuc |E_n|^{1/p} \nrm \chi_{F_n}\nrm_{q^\prime,r^\prime}.
$$
Using this and then H\"{o}lder's inequality in \eqref{localized} (here 
we use  the assumption $q \geq p$), 
we have
\begin{align*}
\langle T\chi_{E }, \chi_{F } \rangle 
& \lesssim \beta^{-\eta} \sum_n \langle T\chi_{E_n}, \chi_{F_n}\rangle 
\\
&\kuc \beta^{-\eta}\sum_n |E_n|^{1/p} \nrm \chi_{F_n}\nrm_{q^\prime,r^\prime}\\
& \kuc \beta^{-\eta}\left[\sum_n |E_n| \right]^{1/p} 
\left[\sum_n   \nrm \chi_{F_n}\nrm_{q^\prime,r^\prime}^{q^\prime}\right]^{1/q^\prime} \\
&\leq \beta^{-\eta} |E|^{1/p} 
    \nrm \chi_{F }\nrm_{q^\prime,r^\prime}.
\end{align*}
In the last inequality we have used the bounded overlap property of the collections 
$\{E_n\}$ and $\{F_n\}$. 
Since $\eta>0$ can be chosen to be arbitrarily small,
this finishes the proof of Theorem~\ref{T:main}.

\section{
Proof of Lemma~\ref{L:balls}}  \label{formerlyAppendix1}
We refer the reader to \cite{christ3} for the proof  of (i) and (ii).

Let $B$ be the ball $B(z_0,\delta_1,\delta_2)$. 
Note that $\{e^{sV_1}z_0:|s|<\delta_1\}\subset B$ and 
by \eqref{v1flow}, $|\Pi\pi_2(\{e^{sV_1}z_0:|s|<\delta_1\})|\approx \delta_1$. 
Now, we prove that 
$|\Pi\pi_2(B)|\kuc \delta_1$. Let $z\in B$ and let $\varphi$ 
be a curve connecting $z_0$ to $z$ as 
in Definition~\ref{D:balls}. Let $\psi=\Pi\circ\pi_2\circ\varphi$. 
Since $V_2$ is contained in the kernel of 
$D\pi_2$, $|a_1(t)|<\delta_1$ and $|V_1|\kuc 1$, we have 
$$
|\psi^\prime(t)| \kuc \delta_1.
$$ 
Therefore $|\Pi\pi_2(B)|\kuc \delta_1$. This proves (iii). 

We have actually shown that the 1-dimensional measure of 
$B(z_0,2\delta_1,2\delta_2) \cap \pi_1^{-1}(x)$ is $\approx \delta_1$ for 
every $x\in \pi_1(B)$. This also proves (ii) for $j=1$. In fact it implies that 
if $A$ is a subset of a $k$-dimensional smooth submanifold of $X$, then
$\delta_1$ times the $k$-dimensional measure of $ A $    is comparable to 
the $k+1$-dimensional measure of $ \pi_1^{-1}(A)\cap B$. The corresponding 
statement holds for $j=2$, as well.

To prove (iv), 
define $f(t,\delta_1,\delta_2) = 
|\pi_2^{-1}\Pi^{-1}(t)\cap B(z_0,\delta_1,\delta_2)|$,
where $|\cdot|$ signifies the $d$-dimensional measure of a subset of $Z$.
To simplify notation we suppose that $\Pi(\pi_2(z_0))=0$; this can be
achieved by a change of coordinates.
Note that if $A\subset B(z_0,\delta_1,\delta_2)$ then for $s < \delta_1$, 
$e^{sV_1}A\subset B(z_0,2\delta_1,2\delta_2)$. 
Choose $t_0\in [-\delta_1, \delta_1]$  so that 
$f(t_0,\delta_1,\delta_2) = \max_{|t|\le\delta_1} f(t,\delta_1,\delta_2)$.
Then
$$
|B(z_0,2\delta_1,2\delta_2)| 
\buy \int_{|s|\le \delta_1} |e^{sV_1}f(t_0)| ds \buy \delta_1 f(t_0).
$$ 
Now $|B(z_0,\delta_1,\delta_2)|\gtrsim |B(z_0,2\delta_1,2\delta_2)|$;
the proof given by Nagel, Stein, and Wainger \cite{nagelsteinwainger}
of the  volume doubling property carries over to two-parameter balls
with $(\theta,A)$--weakly comparable radii, with a bound depending
on $\theta$ and on $A$ but not on $\delta_1,\delta_2,z_0$.
We conclude that
$f(t)\kuc |B(z_0,\delta_1,\delta_2)|/\delta_1$ whenever $|t|\le\delta_1$.

Similar reasoning shows that if $|t|,|t'|\le\delta_1$
then $f(t',\delta_1,\delta_2)\le  C
f(t,3\delta_1,3\delta_2)$. 
Since $\Pi(\pi_2(B(z_0,\delta_1,\delta_2))\subset [-\delta_1,\delta_1]$,
there exists some $t'\in[-\delta_1,\delta_1]$
for which $f(t',\delta_1,\delta_2)\gtrsim 
|B(z_0,\delta_1,\delta_2)|/\delta_1$.
By combining these two observations,
we conclude the reverse inequality
$f(t,3\delta_1,3\delta_2) \gtrsim 
|B(z_0,\delta_1,\delta_2)|/\delta_1$
whenever $|t|\le\delta_1$.

The lower bound for $\|\chi_{\pi_2(B)}\|_{q',r'}$ in (iv) 
follows from this together with the
observation that $f(t)\approx |\pi_2(f(t))| \delta_2$. 
The upper bound follows in the same way from the upper bound $f(t)\lesssim
|B(z_0,\delta_1,\delta_2)|/\delta_1$.

Finally, conclusion (v) follows from (ii) and (iv).
\qed
 
\section{On the hypothesis $p\le q\le r$}   \label{formerlyAppendix2}

One reason why the restriction
$r\geq q \geq p$ in Theorem~\ref{T:main} is natural 
is as follows.

 
\begin{prop} 
Suppose that 
$|B(z,\delta_1,\delta_2)|$
is comparable to 
$|B(z',\delta_1,\delta_2)|$,
uniformly for all $z,z'$ and all weakly comparable $\delta_1,\delta_2$.
Then all valid mixed norm inequalities 
for $T$ are implied by
the conjectured inequalities. That is:\\
(i) If $T$ is of restricted weak mixed type $(p,q,r)$ 
with $r > p >q$, then $(p,p,r)\in P(T)$.\\
(ii)  If $T$ is of restricted weak mixed type $(p,q,r)$ 
with $q> r > p$, then  
$(p,q,r)$ is an interpolant between
$(1,\infty,1)$ and some 
$(p_1,q_1,r_1)\in P(T)$.
\end{prop}

In case (i), the restricted weak mixed type bound
for $(p,p,r)$ implies that for $(p,q,r)$ by H\"older's inequality,
since we are working in a bounded region and $q<p$.
In case (ii), the conclusion is that $(p^{-1},q^{-1},r^{-1})$
belongs to the line segment  with endpoints $(1,0,1)$
and  $(p_1^{-1},q_1^{-1},r_1^{-1})$.
Since any generalized Radon transform $T$ is of strong type $(1,\infty,1)$, 
the restricted weak mixed type $(p,q,r)$ inequality follows from the
$(p_1,q_1,r_1)$ inequality by interpolation.

\begin{proof} (i) If $(p,p,r)\not\in P(T)$ then there exist $\theta,A$ and   a sequence of balls 
$B_n(z_n,\delta_{n,1},\delta_{n,2})$ with $\delta_{n,1}\sim_{(\theta,A)}\delta_{n,2}$ satisfying
$$
|B_n|^{\frac{1}{r}-\frac{1}{p}}\delta_{n,1}^{\frac{2}{p}-\frac{1}{r}}\delta_{n,2}^{1-\frac{1}{r}}\rightarrow \infty,
$$
as $n\rightarrow \infty$. 
Choose
$N_n\approx \delta_{n,1}^{-1}$  balls of size comparable to $B_n$ with   
disjoint projections under $\Pi  \circ  \pi_2$;
this is possible by Lemma~\ref{L:balls}.   
Let $U_n$ be the union of these balls. Note that $|\pi_1(U_n)|\kuc N_n |\pi_1(B_n)|$.
\begin{align*}
\frac{|U_n|}{|\pi_1(U_n)|^{1/p}\nrm\chi_{\pi_2(U_n)}\nrm_{q^\prime,r^\prime}}&\buy
 \frac{N_n|B_n|}{\left[N_n\frac{|B_n| }{\delta_{n,1}}\right]^{1/p}\left[\frac{|B_n|}{\delta_{n,1}\delta_{n,2}}\right]^{1/r^\prime}\delta_{n,1}^{1/q^\prime}N_n^{1/q^\prime}}\\
&\approx \left(|B_n|^{\frac1{r}-\frac1{p}}
\delta_{n,1}^{\frac{2}{p}-\frac{1}{r}}\delta_{n,2}^{1-\frac{1}{r}}\right)
\left(N_n\delta_{n,1}\right)^{1/q-1/p}\\
&\approx |B_n|^{\frac{1}{r}-\frac{1}{p}}\delta_{n,1}^{\frac{2}{p}-\frac{1}{r}}\delta_{n,2}^{1-\frac{1}{r}}\rightarrow \infty,
\end{align*}
as $n\rightarrow \infty$.  Thus, $T$ can not be 
of restricted weak mixed type $(p,q,r)$.

(ii) 
For $s\in [0,1]$ define
$p_1,q_1,r_1$ by
$$
\frac{1}{p }
=1-s+\frac{s}{p_1},\,\,\,\,\,\,\frac{1}{q }=\frac{s}{q_1},\,\,\,\,\,\,
\frac{1}{r}=1-s+\frac{s}{r_1}.
$$
A bit of algebra shows that
$r_1\geq q_1\geq p_1$ if and only if
$\frac{1}{q}+\frac{1}{p^\prime}\leq s \leq \frac{1}{q}+\frac{1}{r^\prime}$,
and moreover that $0<\frac1q + \frac1{p'} \le \frac1q + \frac1{r'}<1$ 
under our assumption that $q>r>p$.
Thus it is possible to choose $s\in(0,1)$ so that
$r_1\geq q_1\geq p_1$. Fix such a parameter $s$.

If $(p_1,q_1,r_1)\in P(T)$ then we have the conclusion of case (ii).
Otherwise
there exist $\theta,A$ and a sequence of  
balls $B_n=B_n(z_n,\delta_{n,1},\delta_{n,2})$ 
with $\delta_{n,1}\sim_{(\theta,A)}\delta_{n,2}$
such that
\begin{align*}
|B_n|^{\frac{1}{r_1}-\frac{1}{p_1}} \delta_{n,1}^{\frac{1}{p_1}+\frac{1}{q_1}-\frac{1}{r_1}} \delta_{n,2}^{1-\frac{1}{r_1}} \rightarrow \infty.
\end{align*}
Note that for these balls
\begin{align*}
\frac{|B_n|}{|\pi_1(B_n)|^{1/p}\nrm\chi_{\pi_2(B_n)}\nrm_{q^\prime,r^\prime}} & \approx
  |B_n|^{\frac{1}{r}-\frac{1}{p}} \delta_{n,1}^{\frac{1}{p}+\frac{1}{q}-\frac{1}{r}} \delta_{n,2}^{1-\frac{1}{r}} 
\\
&= \left(|B_n|^{\frac{1}{r_1}-\frac{1}{p_1}} \delta_{n,1}^{\frac{1}{p_1}+\frac{1}{q_1}-\frac{1}{r_1}} \delta_{n,2}^{1-\frac{1}{r_1}}\right)^s
\rightarrow \infty.
\end{align*}
Thus $T$ can not be of restricted weak mixed type $(p,q,r)$. 
\end{proof}

Case (ii) is valid for all operators $T$,
without the hypothesis that balls of equal bi-radii
have uniformly comparable measures.


\end{document}